\documentclass[10pt,a4paper]{article}
\usepackage{amsfonts}
\usepackage{amssymb}
\usepackage{mathrsfs}
\usepackage{amsthm}
\usepackage{setspace}

\usepackage{fancyhdr}
\newtheorem{defi}{Definition}[section]
\newtheorem{theo}[defi]{Theorem}
\newtheorem{prop}[defi]{Proposition}
\newtheorem{cor}[defi]{Corollary}
\newtheorem{lemm}[defi]{Lemma}

\theoremstyle{definition}
\newtheorem*{rem}{Remark}

\newcommand{\be}{\begin{eqnarray*}}
\newcommand{\ee}{\end{eqnarray*}}
\newcommand{\beqa}{\begin{eqnarray}}
\newcommand{\eeqa}{\end{eqnarray}}
\newcommand{\ba}{\begin{array}}
\newcommand{\ea}{\end{array}}
\newcommand{\rP}{\mathsf{P}}

\newenvironment{lproof}{\emph{Proof of Lemma.}}{ \qed \par}

\begin{document}

\title{Projective Holonomy II: Cones and Complete Classifications}
\author{Stuart Armstrong}
\date{12 September 2006}
\maketitle

\begin{abstract}
The aim of this paper and its prequel is to introduce and classify the irreducible holonomy algebras of the projective Tractor connection. This is achieved through the construction of a `projective cone', a Ricci-flat manifold one dimension higher whose affine holonomy is equal to the Tractor holonomy of the underlying manifold. This paper uses the result to enable the construction of manifolds with each possible holonomy algebra.
\end{abstract}

\section{Introduction}

The aim of this paper is to continue the project started in \cite{mecon}, that of classifying the holonomy algebras of various parabolic geometries. Papers \cite{mecon} and \cite{me!} study conformal holonomies, this one and its predecessor \cite{mepro1} are interested in projective ones. Recall that a projective structure is given by the (unparametrised) geodesics of any affine connection.

In the previous paper \cite{mepro1}, after defining projective structures and introducing the Cartan and Tractor connections, we studied the consequences of reduced projective Tractor holonomy. We found that reducibility on the Tractor bundle $\mathcal{T}$ gave us a foliation of the manifold by Ricci-flat leaves, then showed that the existence of symplectic, complex, hyper-complex and orthogonal structures on the Tractor bundle imply that the underlying manifold is projectively contact, is Einstein, covers a complex space and covers a quaternionic space, respectively. Holonomies of type $\mathfrak{su}$, for instance, correspond to projectively Sasaki-Einstein manifolds. These are not equivalences, however, except in the projectively Einstein case. Table \ref{table:17} gives the list of algebras of projectively Einstein manifolds. Table \ref{table:21} gives the remaining projective holonomy algebras.
\begin{table}[htbp]
\begin{center}
\begin{tabular}{|c|c|c||c|c|}
\hline
\hline
algebra $\mathfrak{g}$ & representation $V$ & \ \ restrictions \ \ & algebra $\mathfrak{g}$ & representation $V$ \\ 
\hline 
\hline
& & & & \\

$\mathfrak{so}(p,q)$ & $\mathbb{R}^{(p,q)}$ & $p+q \geq 5 $ & $\widetilde{\mathfrak{g}}_2$ & $\mathbb{R}^{(4,3)}$ \\ & & & & \\

$\mathfrak{so}(n, \mathbb{C})$ & $\mathbb{C}^{(p,q)}$ & $ n \geq 5 $ & $\mathfrak{g}_2 (\mathbb{C})$ & $\mathbb{C}^7$ \\& & & & \\

$\mathfrak{su}(p,q)$ & $\mathbb{C}^{(p,q)}$ & $p+q \geq 3 $ & $\mathfrak{spin}(7)$ & $\mathbb{R}^8$ \\ & & & & \\

$\mathfrak{sp}(p,q)$ & $\mathbb{H}^{(p,q)}$ & $p+q \geq 2 $ & $\mathfrak{spin}(4,3)$ & $\mathbb{R}^{(4,4)}$ \\ & & & & \\

$\mathfrak{g}_2$ & $\mathbb{R}^7$ & & $\mathfrak{spin}(7, \mathbb{C})$ & $\mathbb{C}^8$ \\ & & & & \\
\hline
\end{tabular}
\end{center}
\caption{Projectively (and Conformally) Einstein Holonomy algebras}
\label{table:17}
\end{table}
\begin{table}[htbp]
\begin{center}
\begin{tabular}{|c|c|c|c|c|}
\hline
\hline
algebra $\mathfrak{g}$ & representation $V$ & \ \ restrictions \ \ & manifold (local) properties \\ 
\hline 
\hline
& & & \\

$\mathfrak{sl}(n, \mathbb{R})$ & $\mathbb{R}^{n}$ & $n \geq 3 $ & Generic \\ & & & \\

$\mathfrak{sl}(n, \mathbb{C})$ & $\mathbb{C}^{n}$ & $ n \geq 3 $ & $U(1)$-bundle over a complex manifold \\& & & \\

$\mathfrak{sl}(n, \mathbb{H})$ & $\mathbb{H}^{n}$ & $ n \geq 2 $ & $Sp(1,\mathbb{H})$-bundle over a quaternionic manifold \\& & & \\

$\mathfrak{sp}(2n,\mathbb{R})$ & $\mathbb{R}^{2n}$ & $n \geq 2 $ & Contact manifold \\ & & & \\

$\mathfrak{sp}(2n,\mathbb{C})$ & $\mathbb{C}^{2n}$ & $n \geq 2 $ & Contact manifold over a complex manifold \\ & & & \\

\hline
\end{tabular}
\end{center}
\caption{Projectively non-Einstein Holonomy algebras}
\label{table:21}
\end{table}

To arrive at these lists, we start by constructing the projective cone: a cone manifold of dimension $n+1$ with an affine connection whose holonomy is isomorphic with that of the Tractor connection of the underlying manifold. This connection is Ricci-flat and torsion-free; thus we may appeal to paper \cite{meric} which, building on \cite{CIH}, gives all possible reductive holonomies for Ricci-flat torsion-free affine connections, and use various tricks and theorems to construct either Ricci-flat cones with the required holonomies, or projective manifolds with the required properties.

These constructions are long and technical, and generate little new mathematics; however they are needed to complete the lists, and a few are interesting in their own rights; the fact that the existence of $\mathfrak{sl}( \mathbb{C})$ type holonomies are much easier to establish than $\mathfrak{sl}(\mathbb{R})$ types is intriguing. Some low-dimensional cases resist the general treatments; these are dealt with individually at the end the paper. Indeed a few Ricci-flat holonomies cannot correspond to Ricci-flat cones at all.

The author would like to thank Dr. Nigel Hitchin, under whose supervision and inspiration this paper was crafted. This paper appears as a section of the author's Thesis \cite{methesis}.

\section{Previous results} \label{old:def}
In the previous paper \cite{mepro1}, we established that the Tractor bundle $\mathcal{T}$ is of rank $n+1$, and a choice of preferred connection -- a torsion-free connection preserving the projective structure -- gives a splitting
\begin{eqnarray*}
\mathcal{T} = T[\mu] \oplus L^{\mu}
\end{eqnarray*}
where $\mu =  \frac{n}{n+1}$, $L^{\mu}$ is the weight bundle $\left( \wedge^n T^* \right)^{\frac{\mu}{-n}}$ and $T[\mu] = T \otimes L^{\mu}$. The preferred connection $\nabla$ acts on $L^{\mu}$ and generates the rho-curvature $\mathsf{P}$ and the Weyl tensor $W$
\begin{eqnarray} \label{pro:rho}
\mathsf{P}_{hj} &=& - \frac{n}{n^2-1} \mathsf{Ric}_{hj} - \frac{1}{n^2 - 1} \mathsf{Ric}_{jh}, \\
\label{pro:weyl}
R_{hj \phantom{k} l}^{\phantom{hj} k} &=& W_{hj \phantom{k} l}^{\phantom{hj} k} + \mathsf{P}_{hl} \delta_j^k + \mathsf{P}_{hj} \delta_l^k - \mathsf{P}_{jl} \delta^k_h - \mathsf{P}_{jh} \delta_l^k.
\end{eqnarray}

The Tractor connection $\overrightarrow{\nabla}$ is given in this splitting by $\overrightarrow{\nabla}_X = \nabla_X + X + \mathsf{P}(X)$, or, more explicitly,
\begin{eqnarray*}
\overrightarrow{\nabla}_X \left( \begin{array}{c} Y \\ a \end{array} \right) = \left( \begin{array}{c}  \nabla_X Y + Xa \\ \nabla_X a + \mathsf{P}(X,Y) \end{array} \right).
\end{eqnarray*}
The curvature of $\overrightarrow{\nabla}$ is then
\begin{eqnarray} \label{curv:def}
R^{\overrightarrow{\nabla}}_{X,Y} = \left( \begin{array}{c} 0 \\ W(X,Y) \\ CY(X,Y) \end{array} \right),
\end{eqnarray}
where $CY$ is the Cotton-York tensor
\begin{eqnarray*}
CY_{hjk} = \nabla_h \mathsf{P}_{jk} - \nabla_j \mathsf{P}_{hk}.
\end{eqnarray*}

\section{Cone construction} \label{cone:construction}
This section will present the main result of this paper: the projective cone construction and its complex analogue. This construction has appeared before in the litterature, such as in Fox's paper \cite{fox} which attributed it to T.~Y.~Thomas \cite{oldcone}, as well as in the author's thesis \cite{methesis}.

First notice that vector line-bundle $L^{-\mu}$ has a principal $\mathbb{R}^+$-bundle -- the quotient of the full frame bundle of $T$ by the action of the simple piece $SL(n, \mathbb{R}) \subset GL(n, \mathbb{R})$. Call this principal bundle $\mathcal{C} (M)$, the cone over $M$. Let $\pi$ be the projection $\pi: \mathcal{C} (M) \to M$ and $Q$ the vector field on $C(M)$ generated by the action of $\mathbb{R}^+$. The main result of this section is:

\begin{theo}
If $(M,\overrightarrow{\nabla})$ is a projective manifold, then there exists a torsion-free Ricci-flat affine connection $\widehat{\nabla}$ on $\mathcal{C}(M)$, which has the same holonomy as $\overrightarrow{\nabla}$.
\end{theo}
This construction bears similarities to the conformal ambient metric construction presented in \cite{ambient1} and \cite{ambient2}; however, instead of using a metric, we shall use the $\mathsf{P}$-tensor, and will not be extending the cone into a second dimension. The rest of this section will be dedicated to proving this.

Fix a preferred connection $\nabla$; this defines not only a splitting of $\mathcal{T}$, but also, because it is a connection on $L^{\mu}$, an $\mathbb{R}^+$-invariant splitting of the projection sequence
\begin{eqnarray*}
0 \longrightarrow \mathbb{R^+} \longrightarrow T \mathcal{C}(M) \stackrel{d \pi}{\longrightarrow} TM \longrightarrow 0.
\end{eqnarray*}
For the rest of this section, let $X$, $Y$ and $Z$ be sections of $TM \subset T\mathcal{C}(M)$. Then define the connection $\widehat{\nabla}$ by
\begin{eqnarray*}
\widehat{\nabla} Q &=& Id, \\
\widehat{\nabla}_Q Y &=& Y, \\
\widehat{\nabla}_{X} Y &=& \nabla_{X} Y + \mathsf{P}(X,Y)Q.
\end{eqnarray*}
\begin{prop}
$\widehat{\nabla}$ is torsion-free, Ricci-flat and projectively invariant.
\end{prop}
\begin{proof}
See \cite{fox} or \cite{methesis}.
\end{proof}

\begin{rem}
One can reconstruct any preferred connections from the cone connection via the following method: given a $Q$-invariant splitting of $T \mathcal{C}(M)$, one has a connection $\nabla$ on $T \subset T \mathcal{C}(M)$ by projecting $\widehat{\nabla}$ along $Q$. In other words
\begin{eqnarray*}
\nabla_X Y = \pi_* (\widehat{\nabla}_X Y).
\end{eqnarray*}
And, of course, $\nabla$ is the preferred connection corresponding to our chosen splitting of $T\mathcal{C}(M)$.
\end{rem}
\begin{rem}
Two such splittings will differ via
\begin{eqnarray*}
X \to X' = X + \Upsilon(X)Q
\end{eqnarray*}
for some one-form $\Upsilon$ on $M$. This is the origin of the fact that two preferred connections differ by the action of a one-form $\Upsilon$.
\end{rem}

\begin{lemm} \label{cone:dif}
Let $\phi$ and $\phi'$ be two paths in $\mathcal{C}(M)$ with identical endpoints such that
\begin{eqnarray*}
\pi (\phi) = \pi(\phi').
\end{eqnarray*}
Then the holonomy transforms of $\widehat{\nabla}$ along $\phi$ and $\phi'$ are the same.
\end{lemm}
\begin{lproof}
Let $X + aQ$ be a vector field, parallel transported along $\phi$,
\begin{eqnarray*}
\widehat{\nabla}_{\dot{\phi}} (X + aQ) = 0.
\end{eqnarray*}
Now there is a (local) invariant extension of $X + aQ$ in the direction of the cone, $e^{-q} (X + aQ)$ where $q$ is a local coordinate, $q = 0$ (locally) along $\phi$ and $Q(q) = 1$. Consequently,
\begin{eqnarray*}
\widehat{\nabla}_Q (e^{-q}(X + aQ)) = 0,
\end{eqnarray*}
and
\begin{eqnarray*}
\widehat{\nabla}_Y (e^{-q}(X + aQ)) = e^{-q} \widehat{\nabla}_Y (X + aQ),
\end{eqnarray*}
so since $\dot{\phi'} = \dot{\phi} + bQ$ for some function $b$,
\begin{eqnarray*}
\widehat{\nabla}_{\dot{\phi'}} (e^{-q}(X + aQ)) &=& e^{-q} \widehat{\nabla}_{\dot{\phi}} (X + aQ) + b\widehat{\nabla}_{Q}(e^{-q}(X + aQ)) \\
&=& 0.
\end{eqnarray*}
Then since $q=0$ locally at both endpoints of $\phi$ and $\phi'$, the result is proved.
\end{lproof}
To complete this section and give a point to it all, one has to show the final result:
\begin{theo}
Differentiating $\mathcal{T}$ along $T$ via $\overrightarrow{\nabla}$ or differentiating $T\mathcal{C}(M)$ along $T \subset T\mathcal{C}(M)$ via $\widehat{\nabla}$ is an isomorphic operation.
\end{theo}
\begin{proof}
A section $s$ of $L^{\mu}$ is isomorphic with a $\mathbb{R}^+$-equivariant function $\mathcal{C}(M) \to \mathbb{R}$. In our case, we require that
\begin{eqnarray*}
Q(s) = \mu s.
\end{eqnarray*}
Then we may identify $(sY, s) \in \Gamma(\mathcal{T})$ with $(sY, sQ) \in \Gamma(T\mathcal{C}(M))$. Under this identification it is clear that
\begin{eqnarray*}
\overrightarrow{\nabla}_X (sY, s) \cong \widehat{\nabla}_X (sY,sQ).
\end{eqnarray*}
\end{proof}
As a simple consequence of this and Lemma \ref{cone:dif},
\begin{cor}
$\overrightarrow{\nabla}$ and $\widehat{\nabla}$ have same holonomy.
\end{cor}
So in order to classify holonomy groups of $\overrightarrow{\nabla}$, one has to look at those groups that can arise as the affine holonomy groups of Ricci-flat cones. By an abuse of notation, so as not to clutter up with too many connection symbols, we will also designate $\widehat{\nabla}$ with the symbol $\overrightarrow{\nabla}$.

Note that if a preferred connection $\nabla$ is Einstein, the projective cone construction is the same as the conformal Einstein cone construction of \cite{mecon}. In that case the conformal $\mathsf{P}^{co}$ is half of the projective $\mathsf{P}$, and the two cone connections are isomorphic.
In this way, by classifying projective Tractor holonomy groups, we shall also classify conformal Tractor holonomy groups for conformally Einstein structures.

\subsection{Complex projective structures} \label{complex:projective}

Let $M^{2n + 1}$ be a projective manifold with a complex structure $J$ on $\mathcal{T}$ -- hence on the cone $\mathcal{C}(M)$. Assume that $\overrightarrow{\nabla}$ is $R$-invariant, where $R = JQ$.
\begin{lemm}
Being $R$-invariant is equivalent to the disappearance of all curvature terms involving $R$.
\end{lemm} \label{R:curve}
\begin{lproof}
$\overrightarrow{R}_{-,-}R = 0$ by definition. Now let $X$ and $Y$ be vector fields commuting with $R$. Then
\begin{eqnarray*}
\overrightarrow{R}_{R,X} Y &=& (\nabla_R \nabla_X - \nabla_X \nabla_R)Y \\
&=& [R, \nabla_X Y].
\end{eqnarray*}
And that expression being zero is precisely what it means for $\nabla$ to be $R$ invariant.
\end{lproof}

Then we may divide out $\mathcal{C}(M)$ by the action of $Q$ and $R$ to get a manifold $N$. Call this projection $\Pi : \mathcal{C}(M) \to N$. Notice that $\Pi$ factors through $M$:
\begin{eqnarray*}
\mathcal{C}(M) \longrightarrow M \longrightarrow N.
\end{eqnarray*}
This makes $\mathcal{C}(M)$ into an ambient construction for the \emph{complex projective structure} detailled in \cite{mepro1}. In brief, $N$ has a well-defined complex strcuture $J_N$, and a host of preferred connections $\widetilde{\nabla}$ that preserve \emph{generalised complex geodesics} on $N$.

We then define the complex cone connection $\mathcal{C}^{\mathbb{C}}(N)$ to be $\mathcal{C}(M)$. The point of this is:
\begin{theo} \label{comproj:proj}
By looking at all possible $(\mathcal{C}(M), \overrightarrow{\nabla}, J)$ that are $R$-invariant, one generates all possible complex projective manifolds $N$. Moreover, $M$ can be reconstructed from $N$.
\end{theo}
\begin{proof}
To prove this, we shall construct a complex cone $\mathcal{C}^{\mathbb{C}}(N)$ for any complex projective manifold $N$. Then $M$ comes directly from dividing $\mathcal{C}^{\mathbb{C}}(N)$ by the action of $Q$.

Given a $N$ with a complex projective structure, choose a preferred connection $\widetilde{\nabla}$; for simplicity's sake, let $\widetilde{\nabla}$ be a preferred connection that preserves a complex volume form. The formulas work for all $\widetilde{\nabla}$, but we won't need that level of generality. Paper \cite{mepro1} defines the complex rho-tensor $\rP^{\mathbb{C}}$.
 
Then let $\mathcal{C}^{\mathbb{C}}(N) = \mathbb{R}^2 \times N$, and let $Q$ and $R$ be the vectors in the direction of $\mathbb{R}^2$. Extend $J^N$ by defining $JQ = R$, and define the connection $\overrightarrow{\nabla}$ as
\begin{eqnarray*}
\overrightarrow{\nabla} Q &=& Id, \\
\overrightarrow{\nabla} R &=& J, \\
\overrightarrow{\nabla}_X Y &=& \widetilde{\nabla}_X Y + \mathsf{P}^{\mathbb{C}}(X,Y)Q - \mathsf{P}^{\mathbb{C}}(X,JY)R,
\end{eqnarray*}
and defining the rest of the terms by torsion-freeness. Then $\overrightarrow{\nabla}$ obviously preserves $J$, and, as in the real projective case,
\begin{lemm}
$\overrightarrow{\nabla}$ is Ricci-flat.
\end{lemm}
\begin{lproof}
Let $X$, $Y$ and $Z$ be sections of $H$. The only non-zero components of the curvature of $\overrightarrow{\nabla}$ is
\begin{eqnarray*}
\overrightarrow{R}_{X,Y}Z &=& \overrightarrow{\nabla}_X \overrightarrow{\nabla}_Y Z- \overrightarrow{\nabla}_Y \overrightarrow{\nabla}_X Z - \overrightarrow{\nabla}_{[X,Y]} Z \\
&=& \widetilde{R}_{X,Y} Z \\
&& + \mathsf{P}^{\mathbb{C}}(Y,Z)X - \mathsf{P}^{\mathbb{C}}(X,Z)Y \\
&& - \mathsf{P}^{\mathbb{C}}(Y,JZ)JX + \mathsf{P}^{\mathbb{C}}(X,JZ)JY  \\
&& + (\nabla_X \mathsf{P}^{\mathbb{C}})(Y,Z)Q - (\nabla_Y \mathsf{P}^{\mathbb{C}})(X,Z)Q \\
&& - (\nabla_X \mathsf{P}^{\mathbb{C}})(Y,JZ)R + (\nabla_Y \mathsf{P}^{\mathbb{C}})(X,JZ)R. \\
\end{eqnarray*}
Most of these terms will disappear upon taking the Ricci trace. In fact
\begin{eqnarray*}
\overrightarrow{\mathsf{Ric}}(X,Y) &=& \widetilde{\mathsf{Ric}}(X,Y) - \mathsf{P}^{\mathbb{C}}(X,Z) + 2n \mathsf{P}^{\mathbb{C}}(X,Z) + \mathsf{P}^{\mathbb{C}}(JX, JZ) \\
&=& \widetilde{\mathsf{Ric}}(X,Z) + (2n-2) \mathsf{P}^{\mathbb{C}}(X,Z) \\
&=& 0.
\end{eqnarray*}
and all other Ricci terms are evidently zero.
\end{lproof}
One may then define the manifold $M$ by dividing out $\mathcal{C}^{\mathbb{C}}(N)$ out by the action of $Q$. Since $\mathcal{C}^{\mathbb{C}}(N)$ is a real cone -- as $\overrightarrow{\nabla} Q = Id$ -- this generates a real projective structure on $M$.
\end{proof}
This real projective structure is $R$-invariant, and generates the original complex projective structure on $N$. This demonstrates that
\begin{prop}
The above construction of $\mathcal{C}^{\mathbb{C}}(N)$ is independent of the choice of $\widetilde{\nabla}$.
\end{prop}
Notice that if a preferred connection $\widetilde{\nabla}$ is holomorphic, then the whole construction is just the complexification of the real case.

\begin{rem}
One may say that a general $\widetilde{\nabla}$ connection gives splittings of $T\mathcal{C}(M)$ and $TM$. If $\widetilde{\nabla}$ preserves a complex volume form up to real multiplication, then the second splitting comes in fact from a section $N \to M$. If $\widetilde{\nabla}$ preserves a complex volume form up to complex multiplication, then the first splitting comes from a section $M \to \mathcal{C}(M)$. And if, as in the example we've dealt with, $\widetilde{\nabla}$ preserves both, then everything is generated by an overall section $N \to \mathcal{C}(M)$.
\end{rem}
\begin{rem}
In terms of splittings of $T\mathcal{C}(M) \cong T\mathcal{C}^{\mathbb{C}}(N)$, a splitting $TN \subset T\mathcal{C}^{\mathbb{C}}(N)$given by a preferred connection $\widetilde{\nabla}$ on $N$ extends to a splitting $TM \subset T\mathcal{C}(M)$ by simply defining
\begin{eqnarray*}
TM = TN \oplus \mathbb{R}(R).
\end{eqnarray*}
The connection corresponding to this splitting is $\nabla$, the $J$-preferred connection that generated $\widetilde{\nabla}$ in the first place \cite{mepro1}.
\end{rem}

Paper \cite{HyQua} details what is actually a quaternionic projective structure, with a hypercomplex cone construction. See the later Section \ref{quat:hol} for more details.

\section{Realisation of holonomy groups} \label{chap:five}
\subsection{Ricci-flat holonomies}
In this section, we will use the list of irreducible holonomies of torsion-free Ricci-flat connections, as established in \cite{meric}, by building on the general torsion-free lists established by \cite{CIH}.
\begin{prop}
Let $\mathfrak{g}$ be a holonomy algebra acting irreducibly on the tangent space. If there exists a torsion-free Ricci-flat connection with holonomy $\mathfrak{g}$, then apart from a few low dimensional exceptions, there exists a projective cone with holonomy $\mathfrak{g}$.
\end{prop}
It is fortunate for our classification result that this is the case, that the holonomy algebra is not an invariant restrictive enough to rule out the cone construction in general.

The list of holonomy algebras permitted by \cite{meric} is:
\begin{eqnarray*}
\begin{array}{|c|c||c|c|}
\hline
\hline
\textrm{algebra }\mathfrak{g} & \textrm{representation V} & \textrm{algebra }\mathfrak{g} & \textrm{representation V} \\
\hline
\hline
& & & \\
\mathfrak{so}(p,q) & \mathbb{R}^{(p,q)}, \ \ p+q \geq 3 & \mathfrak{spin}(3,4)^* & \mathbb{R}^{(4,4)} \\
& & & \\
\mathfrak{so}(n, \mathbb{C}) & \mathbb{C}^n, \ \ n \geq 3 & \mathfrak{spin}(7, \mathbb{C})^* & \mathbb{C}^7 \\
& & & \\
\mathfrak{su}(p,q)^* & \mathbb{C}^{(p,q)}, \ \ p+q \geq 3 & \mathfrak{sl}(n, \mathbb{R}) & \mathbb{R}^n, \ \ n \geq 2 \\
& & & \\
\mathfrak{sp}(p,q)^* & \mathbb{H}^{(p,q)}, \ \ p+q \geq 2 & \mathfrak{sl}(n, \mathbb{C}) & \mathbb{C}^n, \ \ n \geq 1 \\
& & & \\
\mathfrak{g}_2^* & \mathbb{R}^7 & \mathfrak{sl}(n, \mathbb{H})^* & \mathbb{H}^n, \ \ n \geq 1 \\
& & & \\
\widetilde{\mathfrak{g}}_2^{*} & \mathbb{R}^{(3,4)} &\mathfrak{sp}(2n, \mathbb{R}) & \mathbb{R}^{2n}, \ \ n \geq 2 \\
& & & \\
\mathfrak{g}_2(\mathbb{C})^* & \mathbb{C}^7 &\mathfrak{sp}(2n, \mathbb{C}) & \mathbb{C}^{2n}, \ \ n \geq 2  \\
& & & \\
\mathfrak{spin}(7)^* & \mathbb{R}^8 && \\
& & & \\
\hline
\end{array}
\end{eqnarray*}
Algebras whose associated connections \emph{must} be Ricci-flat have been marked with a star.

\begin{rem}
Most constructions in this section will be done by taking the direct product of projective manifolds with known properties. The crux of these ideas is to exploit the fact that projective structures do not respect the taking of direct products: we shall construct examples with maximal Tractor holonomy from the direct product of projectively flat, non-flat manifolds.
\end{rem}

\subsection{Orthogonal holonomy}
The bulk of the work, like the bulk of the possible holonomy groups, lie in this section. We shall construct projective cones for the first ten holonomy algebras.

We will use two approaches: either constructing a projective manifold whose Tractor connection has the holonomy we need, or directly building a projective cone with the required holonomy (and the underlying projective manifold would then emerge by projecting along the cone direction).

\subsubsection{Full orthogonal holonomy}
Here we aim to show that there exist projective manifolds with full $\mathfrak{so}(p,q)$ holonomy algebras. The main theorem is:
\begin{theo} \label{full:holonomy}
Let $(M^m, \nabla^M)$ and $(N^n, \nabla^N)$ be projectively-flat Einstein manifolds, with non-zero Ricci-curvature. Then $(C, \nabla) = (M \times N, \nabla^M \times \nabla^N)$ has full orthogonal holonomy.
\end{theo}
\begin{proof}
Since $M$ is projectively flat, it has vanishing projective Weyl tensor; since it is Einstein, it has symmetric Ricci and rho tensors. Consequently the full curvature of $\nabla^M$ is given by Equation (\ref{pro:weyl}):
\begin{eqnarray*}
(R^M)_{hj \phantom{k}l}^{\phantom{ij}k} = \frac{1}{1-m} \big( \mathsf{Ric}^M_{hl} (\delta^M)^k_j - \mathsf{Ric}^M_{jl} (\delta^M)^k_h \big),
\end{eqnarray*}
with a similar result for $\nabla^N$. Consequently the full curvature of $\nabla$ is
\begin{eqnarray*}
R_{hj \phantom{k}l}^{\phantom{ij}k} = (R^M)_{hj \phantom{k}l}^{\phantom{ij}k} + (R^N)_{hj \phantom{k}l}^{\phantom{ij}k},
\end{eqnarray*}
and its Ricci curvature is
\begin{eqnarray*}
\mathsf{Ric}_{jl}  = \mathsf{Ric}^M_{jl} + \mathsf{Ric}^N_{jl}.
\end{eqnarray*}
Thus the rho tensor of $\nabla$ is
\begin{eqnarray*}
\mathsf{P}_{jl}  = \frac{1}{1-m-n} \big( (1-m) \mathsf{P}^M_{jl} + (1-n) \mathsf{P}^N_{jl}).
\end{eqnarray*}
In other words, the projective Weyl tensor of $(C, \nabla)$ is
\begin{eqnarray} \label{weyl:complicated} \nonumber
&(R^M)_{hj \phantom{k}l}^{\phantom{ij}k} + (R^N)_{hj \phantom{k}l}^{\phantom{ij}k} - \big( \mathsf{P}_{hl} \delta^k_j - \mathsf{P}_{jl} \delta^k_h \big)& \\
&=& \\ \nonumber
& \frac{1}{1-m} \mathsf{Ric}^M_{hl} (\delta^M)^k_j + \frac{1}{1-n} \mathsf{Ric}^N_{hl} (\delta^N)^k_j - \frac{1}{(1-m-n)} \big( \mathsf{Ric}^M_{hl} + \mathsf{Ric}^N_{hl} \big) \delta^k_j &
\end{eqnarray}
minus the corresponding term with $h$ and $j$ commuted. The Cotton-York tensor vanishes, as $\nabla \mathsf{Ric}^M = \nabla \mathsf{Ric}^N = 0$. This expression therefore contains the full curvature of the Tractor connection $\overrightarrow{\nabla}$. Given the splitting defined by $\nabla$,
\begin{eqnarray*}
\mathcal{A} = T^* \oplus \mathfrak{gl}(T) \oplus T,
\end{eqnarray*}
we may start computing the central $(0, \overrightarrow{\mathfrak{hol}}_0, 0) \subset \overrightarrow{\mathfrak{hol}}$ term, by the use of the Ambrose-Singer Theorem \cite{Founddiff} on the Weyl tensor $W$. Because $\nabla$ itself is Einstein (metric $\mathsf{Ric}^M + \mathsf{Ric}^N$, Einstein coefficient one), we know that $\overrightarrow{\mathfrak{hol}}_0 \subset \mathfrak{so}(\mathsf{Ric}^M + \mathsf{Ric}^N)$. Then let
\begin{eqnarray*}
\mu_1 &=& \frac{1}{1-m} - \frac{1}{1-m-n} = \frac{-n}{(1-m)(1-m-n)} \\
\mu_2 &=& \frac{1}{1-n} - \frac{1}{1-m-n} = \frac{-m}{(1-n)(1-m-n)}.
\end{eqnarray*}
Then if $X, Y$ are sections of $T_M$,
\begin{eqnarray*}
W(X,Y) = \mu_1 \mathsf{Ric}^M(X,-)Y - \mu_1 \mathsf{Ric}^M(Y,-)X.
\end{eqnarray*}
Thus $\overrightarrow{\mathfrak{hol}}_0$ must contain $\mathfrak{so}(\mathsf{Ric}^M)$. Similarly for $\mathfrak{so}(\mathsf{Ric}^N)$. These terms lie diagonally inside the maximal bundle:
\begin{eqnarray*}
\left( \begin{array}{cc} \mathfrak{so}(\mathsf{Ric}^M) & 0 \\
0 & \mathfrak{so}(\mathsf{Ric}^M)
\end{array} \right).
\end{eqnarray*}
The upper-right and lower-left components are isomorphic, as representations of $\mathfrak{so}(\mathsf{Ric}^M) \oplus \mathfrak{so}(\mathsf{Ric}^N)$, to $\mathbb{R}^n \otimes \mathbb{R}^m$ and $\mathbb{R}^m \otimes \mathbb{R}^n$, respectively. They are both irreducible as representations, being tensor products of irreducible representations of distinct algebras. Consequently, decomposing $\mathfrak{so}(\mathsf{Ric}^M + \mathsf{Ric}^N)$ as a representation of $\mathfrak{so}(\mathsf{Ric}^M) \oplus \mathfrak{so}(\mathsf{Ric}^N)$, one sees that
\begin{eqnarray*}
\overrightarrow{\mathfrak{hol}}_0 = \mathfrak{so}(\mathsf{Ric}^M) \oplus \mathfrak{so}(\mathsf{Ric}^N) \ \ \textrm{ or } \ \ \overrightarrow{\mathfrak{hol}}_0 = \mathfrak{so}(\mathsf{Ric}^M + \mathsf{Ric}^N),
\end{eqnarray*}
To show that we are in the second case, one merely needs to consider, for $X \in \Gamma(T_M), A \in \Gamma(T_N)$,
\begin{eqnarray*}
W(X,Y) = \frac{-1}{1-m-n} \left( \mathsf{Ric}^M(X,-)A \ - \ \mathsf{Ric}^N(A,-)X \right),
\end{eqnarray*}
evidently not an element of $\mathfrak{so}(\mathsf{Ric}^M) \oplus \mathfrak{so}(\mathsf{Ric}^N)$.

Since $\nabla$ is Einstein, it must preserve a volume form $\nu$, and we know that $\overrightarrow{\nabla}$ preserves a metric $h = \mathsf{Ric}^M + \mathsf{Ric}^N - \nu^2$ on $\mathcal{T}$. The algebra $\mathfrak{so}(h)$ decomposes as $\mathfrak{so}(\mathsf{Ric}^M + \mathsf{Ric}^N) \oplus T$ in terms of the action of $\mathfrak{so}(\mathsf{Ric}^M + \mathsf{Ric}^N)$; the Lie bracket on $\mathfrak{so}(h)$ is given by the natural action of the first component on the latter. Consequently, as before,
\begin{eqnarray*}
\overrightarrow{\mathfrak{hol}} = \mathfrak{so}(\mathsf{Ric}^M + \mathsf{Ric}^N) \ \ \textrm{ or } \ \ \overrightarrow{\mathfrak{hol}} = \mathfrak{so}(h).
\end{eqnarray*}
To show the latter, we turn to infinitesimal holonomy. Since $\nabla$ annihilates both Ricci tensors, we have the expression, for $X, Y, Z$ now sections of $T$:
\begin{eqnarray*}
\overrightarrow{\nabla} (0, W, 0)(X;Y,Z) = (W(Y,Z) \mathsf{P}(X),\  0, \ W(Y,Z)X).
\end{eqnarray*}
And one may evidently choose $X,Y,Z$ to make that last expression non-zero.
\end{proof}

Now we need to find projectively flat manifolds with the required properties. To do so, we define the quadrics
\begin{eqnarray} \label{quad:prop}
S^{(s,t)}(a) = \{ x \in \mathbb{R}^{(s,t)} | g(x,x) = a \}.
\end{eqnarray}
The standard spheres are included in this picture as $S^n = S^{(n+1,0)}(1)$. We may assume $a>0$, as $S^{(s,t)}(a) = S^{(t,s)}(-a)$.

Now $S^{(s,t)}(a)$ is an Einstein manifold with a metric of signature $(s-1,t)$ and positive Einstein coefficient. The $S^{(s,t)}(a)$ are also projectively flat. Using them, we may construct manifolds of dimension $\geq 4$ with orthogonal holonomy of signature $(a,b+1)$ for any non-negative integers $a$ and $b$. However, since orthogonal holonomy with signature $(a,b+1)$ is equivalent to signature $(b+1,a)$, we actually have all the orthogonal holonomy algebras in dimension $\geq 4$.

Consequently
\begin{eqnarray*}
\big( \mathfrak{g}, V \big) \cong \big( \mathfrak{so}(p,q),  \ \mathbb{R}^{(p,q)}, \ p+q \geq 5 \big),
\end{eqnarray*}
are possible projective holonomy algebras.

\begin{theo}
Let $(M^{2m}, \nabla^M)$ and $(N^{2n}, \nabla^N)$ be $\mathbb{C}$-projectively-flat complex Einstein manifolds, with non-zero Ricci-curvature, which is moreover symmetric under the complex structure. Then $(H, \nabla) = (M \times N, \nabla^M \times \nabla^N)$ has full orthogonal $\mathbb{C}$-projective holonomy.
\end{theo}
\begin{proof}
In this case,
\begin{eqnarray*}
\mathsf{P}^{\mathbb{C}}_M &=& \frac{1}{2(1-m)}\mathsf{Ric}^M \\
\mathsf{P}^{\mathbb{C}}_N &=& \frac{1}{2(1-n)}\mathsf{Ric}^N.
\end{eqnarray*}
Now $\nabla= \nabla^N \times \nabla^M$ has Ricci tensor $\mathsf{Ric}^N +\mathsf{Ric}^M$, a symmetric and $\mathbb{C}$-linear tensor; thus $\nabla$ must preserve a complex volume form $\nu$. Then the $\mathbb{C}$-projective holonomy of $M \times N$ must preserve the complex metric $\mathsf{P}^{\mathbb{C}}_M + \mathsf{P}^{\mathbb{C}}_N -\nu^2$. Moreover,
\begin{eqnarray*}
(R^M)_{hj \phantom{k}l}^{\phantom{hj}k} &=& \frac{1}{(1-m)} \big( \mathsf{Ric}^M_{hl} \otimes_{\mathbb{C}} (\delta^M)^k_j - \mathsf{Ric}^M_{jl} \otimes_{\mathbb{C}} (\delta^M)^k_h \big),
\end{eqnarray*}
and similarly for $N$. With these observations, the proof then proceeds in exactly the same way as in the real case.
\end{proof}

To construct such manifolds, one takes the complex versions of the quadrics in the previous argument, and their direct product as before.

By the previous results of Theorem \ref{comproj:proj} any $\mathbb{C}$-projective manifold $M \times N$ with $\mathbb{C}$-projective holonomy algebra $\overrightarrow{\mathfrak{hol}}$ corresponds to a real projective manifold one dimension higher, with $\overrightarrow{\mathfrak{hol}}$ as (real) projective holonomy algebra.

Consequently
\begin{eqnarray*}
\big( \mathfrak{g}, V \big) \cong \big( \mathfrak{so}(n, \mathbb{C}),  \ \mathbb{C}^{n}, \ n \geq 5 \big),
\end{eqnarray*}
are possible projective holonomy algebras.

\subsubsection{$\mathfrak{su}$ holonomies}

When we talk of a manifold with Tractor holonomy $\mathfrak{su}(p,q)$, we are talking about, by definition, a projectively Einstein manifold whose metric cone is Ricci-flat and has holonomy $\mathfrak{su}(p,q)$. In other words this is a Sasaki-Einstein manifold. The existence of Sasaki-Einstein manifolds has been addressed in \cite{SasakiEin} and \cite{SasakiEin2} as well as \cite{suexist}; an adapted proof can also be found in \cite{methesis}, giving all metric signatures needed.

Consequently
\begin{eqnarray*}
\big( \mathfrak{g}, V \big) \cong \big( \mathfrak{su}(m,n),  \ \mathbb{C}^{(m,n)}, \ m+n \geq 3 \big),
\end{eqnarray*}
are possible projective holonomy algebras.

\subsubsection{$\mathfrak{sp}$ holonomies}

When we talk of a manifold with Tractor holonomy $\mathfrak{sp}(p,q)$, we are talking about, by definition, a projectively Einstein manifold whose metric cone is Ricci-flat and has holonomy $\mathfrak{sp}(p,q)$. In other words this is a 3-Sasaki manifold.

The proof of this is similar to the $\mathfrak{su}$ case, except that one uses $N$, an Einstein Quaternionic-K\"ahler, as the base manifold, and $M$ is a principal $SU(2) \cong Sp(1)$ bundle.

Consequently
\begin{eqnarray*}
\big( \mathfrak{g}, V \big) \cong \big( \mathfrak{sp}(m,n),  \ \mathbb{H}^{(m,n)}, \ m+n \geq 3 \big),
\end{eqnarray*}
are possible projective holonomy algebras.

\subsubsection{Exceptional holonomies}
Bryant \cite{MEH} constructs manifolds with exceptional holonomy as cones on other manifolds. All manifolds with exceptional holonomy are Ricci-flat, so these are Ricci-flat cones by definition.

In \cite{MEH}, Bryant shows that the real cone on $SU(3)/T^2$ has holonomy $G_2$ and the real cone on $SU(2,1)/T^2$ has holonomy $\widetilde{G}_2$. Moreover the \emph{complex} cone on $SL(3,\mathbb{C})/T^2_{\mathbb{C}}$ has holonomy $G^{\mathbb{C}}_2$; this corresponds to $SL(3,\mathbb{C})/T^2_{\mathbb{C}}$ having $\mathbb{C}$-projective holonomy $G^{\mathbb{C}}_2$. And, of course, this implies that there exists a manifold one dimension higher -- hence of dimension $15$ -- with real projective holonomy $G^{\mathbb{C}}_2$.

Similarly the cone on $SO(5)/SO(3)$ has holonomy $Spin(7)$. The other $Spin(7)$ cases weren't dealt with in the paper, but one can extend the arguments there to show that the real cone on $SO(3,2)/SO(2,1)$ has holonomy $Spin(3,4)$ and that the complex cone on $SO(5, \mathbb{C})/SO(3, \mathbb{C})$ has holonomy $Spin(7, \mathbb{C})$.

\subsection{Full holonomy}
Here we aim to show that there exist projective manifolds with full $\mathfrak{sl}(n, \mathbb{R})$ holonomy algebras. The main theorem is:
\begin{theo} \label{fullhol:theo}
Let $(M^m, \nabla^M)$ and $(N^n, \nabla^N)$ be projectively-flat manifolds, non-Einstein but with non-degenerate symmetric Ricci tensors. Then $(C = M \times N, \nabla = \nabla^M \times \nabla^N)$ has full holonomy $\mathfrak{sl}(n+m, \mathbb{R})$.
\end{theo}
\begin{proof}
This proof is initially modelled on that of the existence of full orthogonal holonomy in Theorem \ref{full:holonomy}. But first we need:
\begin{lemm}
The Cotton-York tensor of $\nabla$ vanishes.
\end{lemm}
\begin{lproof}
Both manifolds are projectively flat, so have no Tractor curvature. Since the Tractor curvature includes their Cotton-York tensor (Equation (\ref{curv:def})), this last must vanish. So if $X$ and $Y$ are sections of $TM$, $X'$ and $Y'$ sections of $TN$,
\begin{eqnarray*}
(\nabla_X \mathsf{Ric}^M) (Y, -) &=& (\nabla_Y \mathsf{Ric}^M) (X, -) \\
(\nabla_{X'} \mathsf{Ric}^N) (Y', -) &=& (\nabla_{Y'} \mathsf{Ric}^N) (X', -).
\end{eqnarray*}
Then since $\mathsf{Ric}^M$ is covariantly constant in the $N$ direction (and vice versa),
\begin{eqnarray*}
0 = (\nabla_{X'} \mathsf{Ric}^M) (Y, -) &=& (\nabla_Y \mathsf{Ric}^M) (X', -) \\
0 = (\nabla_{X'} \mathsf{Ric}^N) (Y, -) &=& (\nabla_{Y} \mathsf{Ric}^N) (X', -).
\end{eqnarray*}
Consequently the Cotton-York tensor of $(C, \nabla)$ vanishes.
\end{lproof}
Exactly as in Theorem \ref{full:holonomy}, there exists a summand $\mathfrak{h} = \mathfrak{so}(\mathsf{Ric}^M + \mathsf{Ric}^N) \subset \overrightarrow{\mathfrak{hol}}_0$. Under the action of $\mathfrak{h}$, the bundle $\mathfrak{sl}(\mathcal{T})$ splits as
\begin{eqnarray*}
\mathfrak{sl}(\mathcal{T}) = \mathfrak{h} \oplus \odot^2_0 TC \oplus TC \oplus TC^* \oplus \mathbb{R}.
\end{eqnarray*}
Here the bundles $TC$ and $TC^*$ are isomorphic as representations of $\mathfrak{h}$.

Now using infinitesimal holonomy, we consider the first derivative:
\begin{eqnarray} \label{full:realhol}
\overrightarrow{\nabla} \left( \begin{array}{c} 0 \\ W \\ 0 \end{array} \right) (X;Y,Z) &=& \left( \begin{array}{c} W(Y,Z) \mathsf{P}(X) \\  (\nabla_X W)(Y,Z) \\ W(Y,Z)X \end{array} \right).
\end{eqnarray}
Let $X$ and $Y$ be sections of $TM$, $Z$ a section of $TN$. Then Equation (\ref{weyl:complicated}) implies that the central term is
\begin{eqnarray*}
(\nabla_X W)(Y,Z) = \frac{-1}{1-m-n} (\nabla_X \mathsf{Ric}) (Y,-)Z.
\end{eqnarray*}
Since $M$ is non-Einstein, there must exist $X$ and $Y$ such that this term in non-zero. This term is evidently not a section of $\mathfrak{h}$, so
\begin{eqnarray*}
\mathfrak{h} \oplus \odot^2_0 TC = \mathfrak{sl}(TC) \subset \overrightarrow{\mathfrak{hol}}.
\end{eqnarray*}
Now $\mathfrak{sl}(TC)$ does distinguish between $TC$ and $TC^*$; thus looking at Equation (\ref{full:realhol}), we can see that $T \oplus T^* \subset \overrightarrow{\mathfrak{hol}}_0$. Then the last $\mathbb{R}$ term is generated by the Lie bracket between $TC$ and $TC^*$, so
\begin{eqnarray*}
\overrightarrow{\mathfrak{hol}} = \mathfrak{sl}(\mathcal{T}).
\end{eqnarray*}
\end{proof}

We now need to show the existence of such manifolds; in order to do that, we have
\begin{prop}
There exist manifolds with the conditions of Theorem \ref{fullhol:theo}.
\end{prop}
\begin{proof}
Consider $\mathbb{R}^n$, with standard coordinates $x^l$ and frame $X^l = \frac{\partial}{\partial x^l}$ and let $\nabla'$ be the standard flat connection on $N$. Using a one form $\Upsilon$, the connection changes to
\begin{eqnarray*}
\nabla_X Y = \nabla'_X Y + \Upsilon(X)Y + \Upsilon(Y)X.
\end{eqnarray*}
Similarly, since $\nabla'$ is Ricci-flat, the rho-tensor of $\nabla$ is, by paper \cite{mepro1},
\begin{eqnarray*}
\mathsf{P}_{hj} = -\nabla'_{j} \Upsilon_{h} + \frac{1}{2} \Upsilon^2_{hj}.
\end{eqnarray*}
Now if we choose $\Upsilon = dx_1 + \sum_l x_l dx_l$, the tensor $\mathsf{P}$ is given by
\begin{eqnarray*}
\mathsf{P} &=& \sum_l (x_l dx_1 \odot dx_l - dx_l \odot dx_l ) + O(2)
\end{eqnarray*}
This is non-degenerate at the origin. Since $\Upsilon = dx_1 + O(1)$,
\begin{eqnarray*}
\nabla_{X^1} X^2 = X^2 + O(1)
\end{eqnarray*}
and
\begin{eqnarray*}
(\nabla_{X^1} \mathsf{P})(X^2,X^2) &=& X^1 \cdot \mathsf{P}(X^2,X^2) - \mathsf{P}(\nabla_{X^1} X^2,X^2) - \mathsf{P}(X^2, \nabla_{X^1} X^2) \\
&=& 0 - 2 \mathsf{P}(X^2,X^2) \\
&=& -2 + O(1).
\end{eqnarray*}
So $\nabla$ is non-Einstein at the origin. Since being non-degenerate and non-Einstein are open conditions, there exists a neighbourhood of the origin with both these properties. Define this to be $N$. One needs lastly to see that $\mathsf{P}$ (and thus $\mathsf{Ric}$) is symmetric -- equivalently, that $\nabla$ preserves a volume form. One can either see it directly by the formula for $\mathsf{P}$, or one can observe that since $\nabla'$ preserves a volume form, the preferred connection $\nabla$ preserves one if and only if $\Upsilon$ is closed. But this is immediate since
\begin{eqnarray*}
\Upsilon = d \left( x_1 + \sum_l \frac{x_l^2}{2} \right).
\end{eqnarray*}
\end{proof}
Consequently
\begin{eqnarray*}
\big( \mathfrak{g}, V \big) \cong \big( \mathfrak{sl}(n, \mathbb{R}),  \ \mathbb{R}^{n}, \ n \geq 5 \big),
\end{eqnarray*}
are possible projective holonomy algebras.

\subsection{Complex holonomy}
To show that one has full complex holonomy is actually simpler than in the real case. The crucial theorem is:
\begin{theo}
Let $(M^{2m}, \nabla^M)$ and $(N^{2n}, \nabla^N)$ be $\mathbb{C}$-projectively-flat manifolds, both Einstein, with non-degenerate Ricci tensors. Assume further that $\mathsf{Ric}^M$ is $\mathbb{C}$-linear while $\mathsf{Ric}^N$ is $\mathbb{C}$-hermitian. Then $(C = M \times N, \nabla = \nabla^M \times \nabla^N)$ has full complex holonomy $\mathfrak{sl}(n+m, \mathbb{C})$.
\end{theo}
\begin{proof}
In this case,
\begin{eqnarray*}
\mathsf{P}^{\mathbb{C}}_M &=& \frac{1}{2(1-m)}\mathsf{Ric}^M \\
\mathsf{P}^{\mathbb{C}}_N &=& \frac{1}{2(-1-n)}\mathsf{Ric}^N.
\end{eqnarray*}
Consequently the curvature tensors of $\nabla^M$ and $\nabla^N$ are given, according to \cite{mepro1}, by
\begin{eqnarray*}
(R^M)_{hj \phantom{k}l}^{\phantom{hj}k} =& \frac{-1}{m-1}& \big( \mathsf{Ric}^M_{hl} \otimes_{\mathbb{C}} (\delta^M)^k_j - \mathsf{Ric}^M_{jl} \otimes_{\mathbb{C}} (\delta^M)^k_h \big) \\
(R^N)_{hj \phantom{k}l}^{\phantom{hj}k} =& \frac{-1}{n+1} & \big( \mathsf{Ric}^N_{hl} \otimes_{\mathbb{C}} (\delta^M)^k_j - \mathsf{Ric}^N_{jl} \otimes_{\mathbb{C}} (\delta^M)^k_h \\
&& + \mathsf{Ric}^N_{hj} \otimes_{\mathbb{C}} \delta_l^k - \mathsf{Ric}^N_{jh} \otimes_{\mathbb{C}} \delta_l^k \big).
\end{eqnarray*}
As usual, the complex Cotton-York tensor is zero, meaning the full curvature of the Tractor connection is contained in the Weyl tensor. We aim to calculate the $\mathbb{C}$-projective holonomy of $C$. From now on, any implicit tensor product is taken to be complex. Then as in the proof of Theorem \ref{full:holonomy}, it is easy to see that if $X, Y \in \Gamma(TM')$,
\begin{eqnarray*}
W(X,Y) = \mu_1 \mathsf{Ric}^M(X,-)Y - \mu_1 \mathsf{Ric}^M(Y,-)X.
\end{eqnarray*}
Alternatively, if $X, Y \in \Gamma(TN')$,
\begin{eqnarray*}
W(X,Y) &=& \mu_3 \mathsf{Ric}^N(X,-)Y - \mu_3 \mathsf{Ric}^N(Y,-)X \\
&& + \mu_3 \big( \mathsf{Ric}^N(X,Y) - \mathsf{Ric}^N(Y,X) \big),
\end{eqnarray*}
where
\begin{eqnarray*}
\mu_3 = \frac{1}{m+n+1} - \frac{1}{n+1}.
\end{eqnarray*}
Consequently, we can see that
\begin{eqnarray*}
\mathfrak{so}(\mathsf{Ric}^M) \oplus \mathfrak{u}(\mathsf{Ric}^N) \ \subset \ \overrightarrow{\mathfrak{hol}}_0 \ \subset \ \overrightarrow{\mathfrak{hol}}.
\end{eqnarray*}
Where $\overrightarrow{\mathfrak{hol}}_0$ is the $\mathfrak{gl}(TC')$ component of $\overrightarrow{\mathfrak{hol}}$, the $\mathbb{C}$-projective holonomy algebra of $C$. Now under the action of $\mathfrak{so}(\mathsf{Ric}^M) \oplus \mathfrak{u}(\mathsf{Ric}^N)$, $\mathfrak{gl}(TC')$ splits as
\begin{eqnarray*}
\mathfrak{gl}(TC_{\mathbb{C}}) = \left( \begin{array}{c|c} \mathfrak{so}(\mathsf{Ric}^M) \oplus \odot^2_0 TM_{\mathbb{C}} \oplus \mathbb{C} & A \\
\hline B & \mathfrak{u}(\mathsf{Ric}^N) \oplus i \mathfrak{u}(\mathsf{Ric}^N)  \end{array} \right).
\end{eqnarray*}
Here $A = TM' \otimes_{\mathbb{C}} (TN')^*$ and $B = (TM')^* \otimes_{\mathbb{C}} TN^{'}$. These are irreducible, but \emph{not} isomorphic representations of $\mathfrak{so}(\mathsf{Ric}^M) \oplus \mathfrak{u}(\mathsf{Ric}^N)$, because of the action  $\mathfrak{u}(\mathsf{Ric}^N)$. Now if $X$ is a section of $TM'$ and $Y$ is a section of $TN'$,
\begin{eqnarray*}
W(X,Y) = \frac{\mathsf{Ric}^M(X,-)Y}{m+n-1} - \frac{\mathsf{Ric}^N(Y,-)X}{m+n+1},
\end{eqnarray*}
an element of $A \oplus B$ that is neither in $A$ nor in $B$. Consequently $A \oplus B \subset \overrightarrow{\mathfrak{hol}}_0$. But the span of $A \oplus B$ under the Lie bracket is the full algebra $\mathfrak{sl}(TC')$. So
\begin{eqnarray*}
\mathfrak{sl}(TC') \oplus i \mathbb{R} \ \subset \ \overrightarrow{\mathfrak{hol}}_0.
\end{eqnarray*}
Let $\mathfrak{h} = \mathfrak{sl}(TC') \oplus i \mathbb{R}$. Under the action of $\mathfrak{h}$, the full algebra $\mathfrak{sl}(\mathcal{T}, \mathbb{C})$ splits as
\begin{eqnarray*}
\mathfrak{sl}(\mathcal{T}, \mathbb{C}) = \mathfrak{h} \oplus TC' \oplus (TC')^* \oplus \mathbb{R}.
\end{eqnarray*}
\begin{lemm}
$\mathfrak{h} \oplus TC' \oplus (TC')^* \ \subset \ \overrightarrow{\mathfrak{hol}}$
\end{lemm}
\begin{lproof}
This is the standard argument, involving infinitesimal holonomy. $TC'$ and $(TC')^*$ are irreducible non-isomorphic representations of $\mathfrak{h}$. Then
\begin{eqnarray*}
\overrightarrow{\nabla} (0, W, 0)(X;Y,Z) = (W(Y,Z) \mathsf{P^{\mathbb{C}}}(X),\  0, \ W(Y,Z)X).
\end{eqnarray*}
Consequently $TC' \oplus (TC')^* \subset \overrightarrow{\mathfrak{hol}}$.
\end{lproof}
To end the proof, notice that you can generate the final $\mathbb{R}$ term by taking the Lie bracket on $TC' \oplus (TC')^*$. So
\begin{eqnarray*}
\overrightarrow{\mathfrak{hol}} = \mathfrak{sl}(\mathcal{T}, \mathbb{C}).
\end{eqnarray*}
\end{proof}

To construct an explicit example of the previous, it suffices to take $M$ as a complex version of the quadrics of Equation (\ref{quad:prop}), and $N$ to be the (Einstein-K\"ahler) projective plane. As a consequence of this, we have manifolds with full $\mathbb{C}$-projective holonomy, which corresponds, by Theorem \ref{comproj:proj}, to a real projective manifold one dimension higher, with same real projective holonomy algebra.

Consequently,
\begin{eqnarray*}
\big( \mathfrak{g}, V \big) \cong \big( \mathfrak{sl}(n, \mathbb{C}),  \ \mathbb{C}^{n}, \ n \geq 4 \big),
\end{eqnarray*}
are possible projective holonomy algebras.

\subsection{Quaternionic holonomy} \label{quat:hol}
The holonomy algebra $\mathfrak{sl}(n, \mathbb{H})$ forces the manifold to be Ricci-flat by definition \cite{ACHHM} and \cite{quaternionic}, so we shall focus on the cone conditions.

Paper \cite{HyQua}, building on ideas from \cite{DGQM} and \cite{CHQM}, demonstrates that when one has a hypercomplex cone construction $(\mathcal{C}(M), \overrightarrow{\nabla}, I, J, K)$, such that $\overrightarrow{\nabla}$ is invariant under the actions of $IQ$, $JQ$, $KQ$ and -- trivially -- $Q$, one may divide out by these actions to get a Quaternionic manifold $N$. Furthermore, a choice of compatible splitting of $T\mathcal{C}(M)$ is equivalent to a choice of torsion-free connection preserving the quaternionic structure. Thus we have the following natural definitions:
\begin{defi}[Quaternionic Projective Structure]
A quaternionic projective structure is the same as a quaternionic structure on a manifold - a reduction of the structure group of the tangent bundle to $\mathfrak{gl}(n, \mathbb{H}) \oplus \mathfrak{sl}(1, \mathbb{H})$ such that this structure is preserved by a torsion-free connection. The preferred connections are precisely the torsion-free connections preserving this structure. The total space of the cone construction is the bundle
\begin{eqnarray*}
L^{-\frac{n}{n+1}} \otimes H,
\end{eqnarray*}
where $H$ is the natural rank $4$ bundle associated to $\mathfrak{sl}(1, \mathbb{H})$.
\end{defi}
A quaternionic projective structure is thus simply a quaternionic structure. The definition of \cite{quaternionic} for the change of quaternionic connection by a choice of one-form is exactly analogous to our formulas for the change of real or complex preferred connections. See paper \cite{quaternionicweyltensor} for the definition of the quaternionic Weyl tensor (recalling that any quaternionic-K\"ahler manifold is Einstein, so any expression involving the metric can be replaced with one involving the Ricci tensor, for the general case).

In fact, our results are somewhat stronger than in the complex case: since $\overrightarrow{\nabla}$ is hypercomplex,
\begin{eqnarray*}
\overrightarrow{R}_{X,Y} = \overrightarrow{R}_{IX,IY}
\end{eqnarray*}
by \cite{HyQua} and \cite{DGQM}. Consequently all curvature terms involving $IQ$, $JQ$ and $KQ$ vanish and, as in the proof of Lemma \ref{R:curve},
\begin{prop}
Every hypercomplex cone is $IQ$-, $JQ$- and $KQ$-invariant, and thus every hypercomplex cone corresponds to a quaternionic structure.
\end{prop}
Given this definition, one may construct examples similarly to the real and complex cases; indeed:
\begin{theo}
Let $(M^{4m}, \nabla^M)$ and $(N^{4n}, \nabla^N)$ be quaternionicaly-flat manifolds, non-Einstein but with non-degenerate symmetric Ricci tensors. Then $(C = M \times N, \nabla = \nabla^M \times \nabla^N)$ has full quaternionic Tractor holonomy $\mathfrak{sl}(n+m, \mathbb{H})$.
\end{theo}
The proof is analogous to the real case, and one can choose $M$ and $N$ to be quaternionic spaces, with a suitable non-Einstein connection, again as in the real case.

Then one may construct the quaternionic cone and divide out by the action of $Q$ to get a real projective manifold with the same Tractor holonomy.

Consequently,
\begin{eqnarray*}
\big( \mathfrak{g}, V \big) \cong \big( \mathfrak{sl}(n, \mathbb{H}),  \ \mathbb{H}^{n}, \ n \geq 3 \big),
\end{eqnarray*}
are possible projective holonomy algebras.

\subsection{Symplectic holonomies}

The constructions used here were originally discovered in a different context by Simone Gutt, to whom I am very grateful. Paper \cite{Symplectic connections} also contains the construction of what is effectively a `symplectic projective structure', with its own Weyl and rho tensors. Though we will not use or detail this explicitly, it is implicitly underlying some aspects of the present proof.

\subsubsection{Real symplectic} \label{real:symplectic}

Let $M^{2n+1}$ be a contact manifold, with a choice of contact form $\alpha \in \Gamma(T^*)$. We may then define the Reeb vector field $E \in \Gamma(T)$ on $M$ by
\begin{eqnarray*}
\alpha(E) &=& 1, \\
d \alpha (E,-) &=& 0.
\end{eqnarray*}
Since $\alpha$ is a contact form, this suffices to determine $E$ entirely. Let $H \subset T$ be the contact distribution defined by $\alpha (H) =0$.
\begin{lemm}
If $X \in \Gamma (H)$ then $[E, X] \in \Gamma (H)$.
\end{lemm}
\begin{lproof}
By definition,
\begin{eqnarray*}
0 = d \alpha (E, X) &=& E \cdot \alpha(X) - X \cdot \alpha(E) - \frac{1}{2} \alpha([E,X]) \\ &=& -\frac{1}{2} \alpha([E,X]).
\end{eqnarray*}
Hence $[E, X] \in \Gamma (H)$.
\end{lproof}
\begin{lemm}
$\mathcal{L}_E \alpha = 0$.
\end{lemm}
\begin{lproof}
For $X$ a section of $H$,
\begin{eqnarray*}
(\mathcal{L}_E \alpha) (X) &=& \mathcal{L}_E (\alpha(X)) - \alpha(\mathcal{L}_E X) \\
&=& 0 - \alpha([E,X]) = 0.
\end{eqnarray*}
Similarly
\begin{eqnarray*}
(\mathcal{L}_E \alpha) (E) &=& \mathcal{L}_E (\alpha(E)) - \alpha(\mathcal{L}_E E) \\
&=& E \cdot 1 - 0 = 0.
\end{eqnarray*}
\end{lproof}
\begin{lemm} $\mathcal{L}_E (d\alpha) = 0$ \end{lemm}
\begin{lproof}
Immediate since $[ \mathcal{L}, d ] =0$.
\end{lproof}
This gives us the following proposition:
\begin{prop}
Dividing out by the action of the one-parameter sub-group generated by $E$ gives a map $\pi: M \to (N, \nu)$ with $(N, \nu)$ a symplectic manifold and $d \alpha = \pi^* \nu$.
\end{prop}
If $X$ and $Y$ are now sections of $T_N$, they have unique lifts $\overline{X}$ and $\overline{Y}$. Then since $d \alpha (\overline{X}, \overline{Y}) = -\frac{1}{2} \alpha ([\overline{X}, \overline{Y}])$, we have
\begin{eqnarray} \label{torsion:formula}
[\overline{X}, \overline{Y}] = \overline{[X,Y]} -2 \nu(X, Y) E.
\end{eqnarray}
The point of all these constructions in the following theorem:
\begin{theo}
Given $\pi: M \to (N, \nu, \nabla)$ such that $M$ is a contact manifold with contact form $\alpha$ with $d \alpha = \pi^* \nu$, and $\nabla$ a connection preserving $\nu$, there exists a Ricci-flat, torsion-free cone connection $\overrightarrow{\nabla}$ on $\mathfrak{M} = \mathbb{R} \times M$ that preserves the symplectic form $e^{2q} (d \alpha + dq \wedge \alpha)$, where $q$ is the coordinate along $\mathbb{R}$.
\end{theo}
\begin{proof}
Let $s$ be a section of $\odot ^2 TN^*$, $U$ a section of $TN$ and $f$ a function on $N$. Then define the following connection on $\mathfrak{M}$:
\begin{eqnarray} \label{sympcon:prop}
\nonumber \overrightarrow{\nabla}_{\overline{X}} \overline{Y} &=& \overline{\nabla_X Y} \nonumber - \nu(X,Y) E - s(X,Y)Q \\
\nonumber \overrightarrow{\nabla}_E \overline{X} &=& \overrightarrow{\nabla}_{\overline{X}} E = \overline{\sigma X} + \nu(X, U)Q \\
\overrightarrow{\nabla}_E E &=& \pi^* f Q - \overline{U} \\
\nonumber \overrightarrow{\nabla}_Q X &=& X \\
\nonumber \overrightarrow{\nabla}_Q E &=& E \\
\nonumber \overrightarrow{\nabla} Q &=& Id.
\end{eqnarray}
Where $s (X,Y)= \nu (X, \sigma Y)$, or, in other words, $\sigma_j^k =s_{hj} \nu^{hk}$. One can see immediately from Equation (\ref{torsion:formula}) that $\overrightarrow{\nabla}$ is torsion-free. It is obviously a cone connection. On top of that:
\begin{prop}
$\overrightarrow{\nabla}$ is a symplectic connection, for the non-degenerate symplectic form $\omega = e^{2q} (d \alpha + dq \wedge \alpha)$.
\end{prop}
\begin{proof}
By direct calculation.
\end{proof}
We may now calculate the curvature of $\overrightarrow{\nabla}$; it is, for $R$ the curvature of $\nabla$,
\begin{eqnarray}
\nonumber \overrightarrow{R}_{Q,-} - &=& \overrightarrow{R}_{-,-} Q = 0 \\
\nonumber \overrightarrow{R}_{X,Y} Z &=& R_{X,Y} Z + 2\nu(X,Y) \sigma Z + \\
\nonumber && (( \sigma Y) \nu(X,Z) - (\sigma X) \nu(Y,Z) + Y s(X,Z) - X s(Y,Z)) \\
\nonumber  && + \left( \nu(X,D(Y,Z)) -\nu(Y,D(X,Z)) \right) Q \\
\label{symp:curvexpr} \overrightarrow{R}_{X,Y} E &=& D(X,Y) - D(Y,X) \\
\nonumber && + (\nu(Y, \nabla_X U) - \nu(X,\nabla_Y U) +2f \nu(X,Y)) Q \\
\nonumber \overrightarrow{R}_{X,E} Y &=& D(X,Y) + (\nu(Y, \sigma^2 X) + f \nu(X,Y) + \nu(Y,\nabla_X U))Q \\
\nonumber \overrightarrow{R}_{X,E} E &=& ( X.f + s(X,U) - \nu(\sigma X,U)) Q\\
\nonumber  && +fX -\nabla_X U -\sigma^2X.
\end{eqnarray}

Where $D(X,Y) = (\nabla_X \sigma )(Y) + \nu(Y,U)X - \nu (X,Y)U$. Taking traces, with $\mathsf{Ric}$ the Ricci curvature of $\nabla$,
\begin{eqnarray*}
\overrightarrow{\mathsf{Ric}}_{-,Q} &=& \overrightarrow{\mathsf{Ric}}_{Q,-} = 0 \\
\overrightarrow{\mathsf{Ric}}_{X,Z} &=& \mathsf{Ric}_{X,Z} + \\ && (\mathrm{trace} \ \sigma ) \nu(X,Z) + 3\nu(\sigma X, Z) +(1-2n) s(X,Z) \\
\overrightarrow{\mathsf{Ric}}_{X,E} &=& -\mathrm{trace} \ (\nabla_X \sigma) + i(X). \mathrm{trace} [ Y \to \nabla_Y \sigma] + (2n+1) \nu(X,U) \\
\overrightarrow{\mathsf{Ric}}_{E,X} &=& i(X). \mathrm{trace} [ Y \to \nabla_Y \sigma] + (2n+1) \nu (X,U) \\
\overrightarrow{\mathsf{Ric}}_{E,E} &=& f 2n - \mathrm{trace} \ (\nabla U) - \mathrm{trace} \ (\sigma^2).
\end{eqnarray*}

Now $\nu( \sigma X, Y) = \nu_{km} \sigma^k_j X^j Y^m = \nu_{km} s_{ij}\nu^{ik} X^j Y^m = -s(X,Y)$, and trace $\sigma =$ trace $\nabla_X \sigma = 0$, so $\overrightarrow{\mathsf{Ric}}$ is symmetric, as expected. So the full expression is:
\begin{eqnarray*}
\overrightarrow{\mathsf{Ric}}_{-,Q} &=& \overrightarrow{\mathsf{Ric}}_{Q,-} = 0 \\
\overrightarrow{\mathsf{Ric}}_{X,Z} &=& \mathsf{Ric}_{X,Z} - (2n+2) s(X,Z) \\
\overrightarrow{\mathsf{Ric}}_{X,E} &=& + i(X). \mathrm{trace} [ Y \to \nabla_Y \sigma] + (2n+1) \nu(X,U) \\
\overrightarrow{\mathsf{Ric}}_{E,X} &=& +  i(X). \mathrm{trace} [ Y \to \nabla_Y \sigma] + (2n+1) \nu (X,U) \\
\overrightarrow{\mathsf{Ric}}_{E,E} &=& f 2n - \mathrm{trace} \ (\nabla U) - \mathrm{trace} \ (\sigma^2).
\end{eqnarray*}
Choose $s = \frac{1}{2n+2} \mathsf{Ric}$, and, for $\eta = \mathrm{trace} [ Y \to \nabla_Y \sigma]$, define $U = - \frac{1}{2n+1} \nu(\eta,-)$. Finally, let $f = \frac{1}{2n} (\mathrm{trace} \ (\nabla U) + \mathrm{trace} \ (\sigma^2))$. Then $\overrightarrow{\nabla}$ is Ricci-flat, as theorised.
\end{proof}

\begin{rem}
It is nearly certainly not the case, however, that \emph{every} Ricci-flat symplectic cone connection can be generated in the above manner; for the $\overrightarrow{\nabla}$ so generated is $E$ invariant, which is not a general condition for a symplectic connection.
\end{rem}
We now aim to construct an explicit connection $\nabla$ such that the $\overrightarrow{\nabla}$ it generates has maximal holonomy.

Let $V$ be the standard representation of $\mathfrak{g} = \mathfrak{sp}(2n, \mathbb{R})$. Then $\mathfrak{g}$ is isomorphic, via the alternating form $\nu$, with $V^* \odot V^*$. The Lie bracket is given, in terms of this isomorphism, as
\begin{eqnarray*}
[ ab, cd ] = \nu (b,c) a d + \nu(b,d)ac + \nu(a,c) bd + \nu(a,d)bc.
\end{eqnarray*}
We know that $H^{(1,2)}(\mathfrak{g}) = 0$ and that all symplectic structures are flat. Moreover $\mathfrak{g}^{(1)} =  \odot^3 V^*$, means that any symplectic connection is locally isomorphic with a section $U \to \odot^3 TN^*$, for $U \subset N$ open. Choosing local symplectic coordinates $(x_j)$ such that
\begin{eqnarray*}
\nu = dx_1 \wedge dx_2 + dx_3 \wedge dx_4 + \ldots + dx_{2n-1} \wedge dx_{2n},
\end{eqnarray*}
we may define the symplectic connection $\nabla$ as
\begin{eqnarray*}
\nabla = d + \sum_{j \neq 1} x_1 dx_1 (dx_j)^2 + \sum_{k \neq 1,2} x_2 dx_2 (dx_k)^2.
\end{eqnarray*}
Notice that $\nabla = d + O(1)$. We may calculate the curvature of $\nabla$ as
\begin{eqnarray*}
R &=& \nabla \wedge \nabla \\
&=& 2 \sum_{j \neq 1} (dx_1 \wedge dx_j) \otimes (dx_1 dx_j) + 2 \sum_{k \neq 1,2} (dx_2 \wedge dx_k) \otimes (dx_2 dx_k) + O(2).
\end{eqnarray*}
When taking the Ricci trace using the symplectic form $\nu$, all terms apart from $(dx^1 \wedge dx^2) \otimes (dx^1 dx^2)$ vanish. Consequently the Ricci tensor is
\begin{eqnarray*}
\mathsf{Ric} = 2 (dx_1 dx_2) + O(2).
\end{eqnarray*}
And, of course,
\begin{eqnarray*}
\nabla \mathsf{Ric} = O(1).
\end{eqnarray*}
This allows us to simplify the curvature equations. By definition $U = O(1)$, so
\begin{eqnarray*}
\overrightarrow{R}_{X,Y} Z &=& R_{X,Y} Z + 2\nu(X,Y) \sigma Z + O(1)\\
&& (( \sigma Y) \nu(X,Z) - (\sigma X) \nu(Y,Z) + Y s(X,Z) - X s(Y,Z)) \end{eqnarray*}
$s = \frac{1}{2n+2}\mathsf{Ric}$ as before.
\begin{prop}
$\overrightarrow{\nabla}$ has full symplectic holonomy.
\end{prop}
\begin{proof}
Still working in our chosen basis, we notice that because of our conditions on the Ricci tensor, for one of $j$ and $k$ in the set $(1,2)$ but $(j,k) \neq (1,2)$,
\begin{eqnarray*}
\overrightarrow{R}_{X^j,X^k} Z = R_{X^j,X^k} Z + O(1)
\end{eqnarray*}
where $X^j = \frac{\partial}{\partial x^j}$. This means, by the Ambrose-Singer Theorem \cite{Founddiff}, that elements of the form $dx_1 dx_j |_0$, $j \neq 2$ and $dx_2 dx_k |_0$, $k \neq 1$, are contained in $\overrightarrow{\mathfrak{hol}}$, the infinitesimal holonomy algebra of $\overrightarrow{\nabla}$ at $0$. Now we may take a few Lie brackets:
\begin{eqnarray} \label{bracket:symp}
[ dx_1 dx_k, dx_2 dx_j] = dx_k dx_j + \nu(dx_k, dx_j)dx_1 dx_2 + O(1)
\end{eqnarray}
implying
\begin{eqnarray*}
[ dx_1 dx_k, dx_2 dx_j] - [dx_1 dx_j, dx_2 dx_k] = 2 dx_k dx_j + O(1).
\end{eqnarray*}
Consequently $dx_k dx_j |_0 \in \overrightarrow{\mathfrak{hol}}$. By (\ref{bracket:symp}), we also have $dx_1 dx_2 |_0$ in this bundle. To show that we have all of $\mathfrak{sp}(\nu, \mathbb{R})$, we need only to add the elements $dx_2 dx_2 |_0$ and $dx_1 dx_1 |_0$. These are generated, for $j$ odd, by
\begin{eqnarray*}
\big[ dx_1 dx_{j}, dx_1 dx_{j+1} \big] &=& dx_1 dx_1 + O(1) \\
\big[ dx_2 dx_{j}, dx_2 dx_{j+1} \big] &=& dx_2 dx_2 + O(1).
\end{eqnarray*}
Under the action of $\mathfrak{sp}(\nu, \mathbb{R})$, the full algebra
$\mathfrak{sp}(\omega, \mathbb{R})$ splits as
\begin{eqnarray*}
\mathfrak{sp}(\omega, \mathbb{R}) = \mathfrak{sp}(\nu, \mathbb{R}) \oplus 2 V \oplus \mathfrak{sp}(\omega/\nu, \mathbb{R}),
\end{eqnarray*}
where the last module is a trivial representation for $\mathfrak{sp}(\nu, \mathbb{R})$.
\begin{lemm}
If $\overrightarrow{\mathfrak{hol}}$ acts irreducibly on $T\mathfrak{M}|_0$, then $\overrightarrow{\mathfrak{hol}} = \mathfrak{sp}(\omega, \mathbb{R})$.
\end{lemm}
\begin{lproof}
If $\overrightarrow{\mathfrak{hol}}$ acts irreducibly on $T\mathfrak{M}|_0$, then
\begin{eqnarray*}
\mathfrak{sp}(\nu, \mathbb{R}) \oplus 2 V \subset \overrightarrow{\mathfrak{hol}},
\end{eqnarray*}
and the $2V$ generate the remaining piece $\mathfrak{sp}(\omega/\nu, \mathbb{R})$ through the Lie bracket.
\end{lproof}
So in order to finish this proof, we need to show that $\overrightarrow{\mathfrak{hol}}$ acts irreducibly on $T\mathfrak{M}_0$, or equivalently,
\begin{lemm}
If there exists $\mathcal{K} \subset T\mathfrak{M}$ with $TN_0 \subset \mathcal{K}_0$ and such that $\mathcal{K}$ is preserved by $\overrightarrow{\nabla}$, then $\mathcal{K} = T\mathfrak{M}$.
\end{lemm}
\begin{lproof}
First of all $\mathcal{K}$ has a non-trivial intersection with $TN$ away from $0$. So let $s \in \Gamma(K) \cap TN$ such that $s(0) = X^1$. Then by Equation (\ref{sympcon:prop})
\begin{eqnarray*}
(\overrightarrow{\nabla}_{X^1} s)(0) = \frac{1}{2n+1} Q + t,
\end{eqnarray*}
whereas
\begin{eqnarray*}
(\overrightarrow{\nabla}_{X^2} s)(0) = E + t',
\end{eqnarray*}
where $t, t' \in TN_0 \subset \mathcal{K}_0$. Consequently $\mathcal{K} = T\mathfrak{M}$, and the holonomy algebra of $\overrightarrow{\nabla}$ acts irreducibly on $T \mathfrak{M}$.
\end{lproof}
\end{proof}

Consequently
\begin{eqnarray*}
\big( \mathfrak{g}, V \big) \cong \big( \mathfrak{sp}(2n, \mathbb{R}),  \ \mathbb{R}^{2n}, \ n \geq 3 \big),
\end{eqnarray*}
are possible projective holonomy algebras.

\subsubsection{Complex symplectic}
The previous proof works exactly the same in the holomorphic category.

Consequently
\begin{eqnarray*}
\big( \mathfrak{g}, V \big) \cong \big( \mathfrak{sp}(2n, \mathbb{C}),  \ \mathbb{C}^{2n}, \ n \geq 3 \big),
\end{eqnarray*}
are possible projective holonomy algebras.

\subsection{Low-dimension cases}

Some low-dimensional algebras are possible affine holonomy algebras, but have not yet been either constructed or ruled out as normal Tractor holonomy algebras. They are:
\begin{eqnarray*}
\begin{array}{|c|c|c|}
\hline
\hline
\textrm{algebra }\mathfrak{g} & \textrm{representation V} & \textrm{Dimensions} \\
\hline
\hline
& & \\
\mathfrak{so}(p,q) & \mathbb{R}^{(p,q)} &  p+q = 3,4  \\
\mathfrak{so}(n, \mathbb{C}) & \mathbb{C}^n & n = 3,4 \\
\mathfrak{sp}(p,q) & \mathbb{H}^{(p,q)} & p+q = 2^* \\
\mathfrak{sl}(n, \mathbb{R}) & \mathbb{R}^n & n= 2,3^*,4^* \\
\mathfrak{sl}(n, \mathbb{C}) & \mathbb{C}^n & n=  1,2,3^* \\
\mathfrak{sl}(n, \mathbb{H}) & \mathbb{H}^n & n=1,2^* \\
\mathfrak{sp}(2n, \mathbb{R}) & \mathbb{R}^{2n} & n = 2^* \\
\mathfrak{sp}(2n, \mathbb{C}) & \mathbb{C}^{2n} & n = 2^*  \\
& & \\
\hline
\end{array}
\end{eqnarray*}
Those marked with stars are those algebras that \emph{can} appear as projective normal Tractor holonomy algebras.

\begin{prop}
The low-dimensional $\mathfrak{so}$ algebras cannot appear as projective holonomy algebras.
\end{prop}
\begin{proof}
Dimensional considerations imply that the conformal Weyl tensor vanishes in $3$ dimensions, \cite{WT}. The obstruction to conformal flatness is carried entirely by the Cotton-York tensor, which of course vanishes for an Einstein space.

So any $3$-dimensional Einstein space is conformally flat -- hence projectively flat, since the two cones are the same. This eliminates the real $\mathfrak{so}$ and the $\mathfrak{su}$, as the underlying manifold must be Einstein. The complex $\mathfrak{so}$ has $\mathbb{C}$-linear curvature, so is automatically holomorphic -- so disappears just as in the real case, as the holomorphic Weyl tensor must also vanish in three complex dimensions.

In two dimensions, Cartan connections no longer correspond to conformal structures, but rather to M\"obius structures \cite{Mob}. The projective cone construction for an Einstein connection is then equivalent to a M\"obius structure which preserves a Tractor -- it is not hard to see that this is flat, see for example \cite{methesis}.

\end{proof}
Since every one-dimensional manifold is projectively flat, $\mathfrak{sl}(2, \mathbb{R})$ and $\mathfrak{sl}(1, \mathbb{C})$ are not possible Tractor holonomy algebras -- they are not even possible Ricci-flat algebras, in fact.
\begin{lemm} \label{sl:two}
$\mathfrak{sl}(2, \mathbb{C})$ is not a possible Tractor holonomy algebra.
\end{lemm}
\begin{lproof}
Assume that $\overrightarrow{\nabla}$ is a cone connection with this holonomy, and let $R = JQ$. From Lemma \ref{R:curve}, we know that a cone connection is $R$-invariant if and only if all curvature terms involving $R$ vanish. For a connection with holonomy $\mathfrak{sl}(2, \mathbb{C})$, being Ricci-flat is equivalent to having $J$-hermitian curvature. Consequently
\begin{eqnarray*}
\overrightarrow{R}_{R,X} = \overrightarrow{R}_{JR,JX} = -\overrightarrow{R}_{Q,JX} = 0.
\end{eqnarray*}
So $\overrightarrow{\nabla}$ is $R$-invariant, and, as in Section \ref{complex:projective}, there is a complex projective manifold $N$ of complex dimension one, for which $\overrightarrow{\nabla}$ is the complex cone connection.

Any two torsion-free complex connections $\widetilde{\nabla}$ and $\widetilde{\nabla}'$ on $N$ differ by a one-form $\Xi \in \Omega^{(1,0)}(N)$
\begin{eqnarray*}
\widetilde{\nabla}_X Y = \widetilde{\nabla}'_X Y + \Xi(X) Y
\end{eqnarray*}
Since $\Xi(X) Y = \Xi(Y) X$, we can set $\Upsilon^{\mathbb{C}} = \frac{1}{2} \Xi$ to see that $\widetilde{\nabla}$ and $\widetilde{\nabla}'$ define the same complex projective structure. So every complex projective structure on $N$ is flat, implying that $\overrightarrow{\nabla}$ itself must be flat.
\end{lproof}
\begin{lemm}
$\mathfrak{sl}(1, \mathbb{H})$ is not a possible Tractor holonomy algebra.
\end{lemm}
\begin{lproof}
We know what $\overrightarrow{\nabla}$ must be, explicitly; it is given by one vector field $Q$ with $\overrightarrow{\nabla} Q = Id$, and three (non-commuting) vector fields $J_{\alpha}Q$ such that
\begin{eqnarray*}
\overrightarrow{\nabla} J_{\alpha}Q = J_{\alpha}.
\end{eqnarray*}
and
\begin{eqnarray*}
[ J_{\alpha}Q, J_{\beta}Q] = -2 J_{\alpha} J_{\beta} Q
\end{eqnarray*}
whenever $\alpha \neq \beta$. But in this case all the curvature of $\overrightarrow{\nabla}$ vanishes.
\end{lproof}
All other low-dimension algebras are possible Tractor holonomies:
\begin{prop}
The algebras $\mathfrak{sl}(2, \mathbb{H})$ and $\mathfrak{sp}(p,q)$, $p+q = 2$, do exist as Tractor holonomy algebras.
\end{prop}
\begin{proof}
As seen in Section \ref{quat:hol}, any hypercomplex Tractor connection corresponds to a quaternionic structure on a manifold $N$, in other words a $\mathfrak{g} = \mathbb{R} \oplus \mathfrak{sl}(1, \mathbb{H}) \oplus \mathfrak{sl}(1, \mathbb{H})$ structure. However this last algebra is equal to $\mathfrak{co}(4)$ -- if one takes $\mathbb{H}$ as a model space, a definite-signature metric $g$ is given by $g(a,b) = Re(a \overline{b})$, and it is easy to see that $\mathfrak{g}$ preserves $g$ up to scaling.

As usual, a subgroup of $\mathfrak{sl}(2, \mathbb{H})$ acting reducibly on $\mathbb{H}^2$ corresponds to a conformally Ricci-flat $4$-fold. But, from results on conformal holonomy in paper \cite{mecon}, we know there exist non conformally Ricci-flat manifolds in four dimensions. The subgroups acting irreducibly on $\mathbb{H}^2$ are
\begin{eqnarray*}
\mathfrak{sp}(2,0) \cong \mathfrak{sp}(0,2), \ \mathfrak{sp}(1,1), \ \mathfrak{sl}(2, \mathbb{H}),
\end{eqnarray*}
corresponding respectively to conformally Einstein with $\lambda <0$, conformally Einstein with $\lambda >0$, and not conformally Einstein at all. Examples of all these constructions, without further holonomy reductions, exist in four dimensions, see Theorem \ref{full:holonomy} and Equation \ref{quad:prop}.
\end{proof}
\begin{rem}
The argument for the rest of this section can be paraphrased as `if we have a manifold with non-trivial Tractor holonomy, we can conjugate the holonomy algebra by gluing the manifold to a copy of itself with a twist, to generate the full algebra'. The subtleties will be in making the manifold flat around the gluing point. This argument only works if the flattening respects whatever structures -- complex or symplectic -- we are attempting to preserve. We must also avoid using Ricci-flat connections, as then conjugation will not give us the full algebras; but it is simple to pick a preferred connection that is not Ricci-flat.
\end{rem}

\begin{prop} \label{sl:three}
The algebras $\mathfrak{sl}(3, \mathbb{R})$ and $\mathfrak{sl}(3, \mathbb{C})$, do exist as Tractor holonomy algebras.
\end{prop}
\begin{proof}
The projective Weyl tensor vanishes in two real dimensions, and consequently the full obstruction to projective flatness is carried by the Cotton-York tensor (see Equation (\ref{curv:def})). Cartan \cite{ECP} proved propositions about two dimensional projective structures that are equivalent to stating that the only possible tractor holonomy algebras are $\mathfrak{sl}(3, \mathbb{R})$ and $\mathfrak{sl}(2, \mathbb{R}) \rtimes \mathbb{R}^{2*}$.

In order to prove the existence of a manifold with full tractor holonomy, we shall use the following proposition:
\begin{prop} \label{patch:U}
Assume there exists a manifold $M^2$ with non-trivial Tractor holonomy. Then there exists a manifold $N^2$ with full Tractor holonomy.
\end{prop}
\begin{proof}
As we've seen, the Tractor holonomy of $M$ is $\mathfrak{sl}(3, \mathbb{R})$ or $\mathfrak{h} = \mathfrak{sl}(2, \mathbb{R}) \rtimes \mathbb{R}^{2*}$. Since the former gives us our result directly, assume the latter; since a non-trivial holonomy algebra must be non-trivial on some set, we have a set $U \subset M$ such that the local holonomy at any point of $U$ is $\mathfrak{h}$. Let $\nabla$ be any preferred connection of this projective structure. Choose local coordinates $(x_j)$ on $U \subset M$, and let $\widehat{\nabla}$ be the flat connection according to these local coordinates. Let $f$ be a bump function, and define $\nabla' = f \nabla + (1-f) \widehat{\nabla}$. This is a torsion free connection, and since $\nabla' =  \nabla$ where $f = 1$, has Tractor holonomy containing $\mathfrak{h}$.

Take two copies $U_1$ and $U_2$ of $(U,\nabla', x_j)$ and identify two small flat patches of them -- patches where $\nabla'_1$ and $\nabla'_2$ are flat --  using the rule $x_j - a_j \to s(x_j - b_j)$ for $s$ some element of $SL(2,\mathbb{R})$, $a$ a point in the flat part of $U_1$, $b$ a point in the flat part of $U_2$, and the $x_j$ local, flat coordinates. This identifies flat sections with flat sections, so does not affect the local holonomy around these patches. The local derivative of $s$ is $Ds(X^j) = s(X^j)$.

Restrict $U_1$ and $U_2$ so that the construction we get is a manifold. Since $s$ maps flat sections to flat sections, $\nabla'_1 = \nabla'_2$ whenever they are both defined. So we have a globally defined $\nabla'$.

Changing $s$ changes the inclusion of the holonomy-preserved vector from $U_1$ into $U_2$, thus changes the inclusion $\mathfrak{h} \subset \overrightarrow{\mathfrak{hol}}_x$ by conjugation on the $\mathfrak{sl}(2, \mathbb{R})$ factor of $\mathcal{A}_x$ defined by $\nabla'$. But any two conjugate non-identical copies of $\mathfrak{h}$ generate all of $\mathfrak{sl}(3,\mathbb{R})$, so we are done.
\end{proof}
Then we may conclude with the following lemma:
\begin{lemm} \label{exist:U}
There exists manifolds $M^2$ with non trivial Tractor holonomy.
\end{lemm}
\begin{lproof}
To do so, it suffices to find a manifold with non-trivial Cotton-York tensor. But if we have local coordinates $x$ and $y$ and corresponding vector fields $X$ and $Y$. Define $\nabla$ such that $\nabla_X Y = \nabla_Y X = \nabla_Y Y = 0$ and $\nabla_X X = y^2 Y$. $\nabla$ is torsion-free and
\begin{eqnarray*}
\mathsf{Ric}^{\nabla} = 2y \ dx \otimes dx,
\end{eqnarray*}
thus
\begin{eqnarray*}
CY^{\nabla} = 4 \ dx \wedge dy \otimes dx.
\end{eqnarray*}
\end{lproof}
We may use these same ideas to construct a manifold with complex projective Tractor holonomy $\mathfrak{sl}(3, \mathbb{C})$ -- and hence a real projective manifold with same holonomy, one dimension higher. The existence proof Lemma \ref{exist:U} works in the holomorphic category, and in then has a tractor holonomy algebra containing $h \otimes \mathbb{C}$. Then given a holomorphic $M$ with these properties, we can use the trick of Proposition \ref{patch:U}, with $(x_j)$ holomorphic coordinates, to get $\nabla' = f \nabla + (1-f) \widehat{\nabla}$. This obviously preserves the complex structure (though it is not holomorphic), and we can then patch $U_1^{\mathbb{C}}$ and $U_2^{\mathbb{C}}$ together using $s \in SL(2, \mathbb{C})$, which also preserves the complex structure.
\end{proof}
\begin{cor}
Tractor holonomy $\mathfrak{sl}(4, \mathbb{C})$ also exists.
\end{cor}
\begin{proof}
The cone over any manifold with Tractor holonomy $\mathfrak{sl}(3,\mathbb{R})$ has Tractor holonomy $\mathfrak{sl}(3,\mathbb{R})$ as well (to see this, see the properties of projectively Ricci-flat manifolds in paper \cite{mepro1}, which demonstrate the cone has Tractor holonomy $\mathfrak{sl}(3,\mathbb{R})$ or $\mathfrak{sl}(3,\mathbb{R}) \rtimes \mathbb{R}^3$; the properties of a cone insure the former. The author's thesis \cite{methesis} shows this in detail). Then we construct the cone over the manifold of the previous proposition, choose a preferred connection that does \emph{not} make it Ricci flat (so that the tangent bundle $T[\mu]$ is not holonomy preserved), and then use the same patching process to conjugate $\mathfrak{sl}(3,\mathbb{R})$ and get full Tractor holonomy.
\end{proof}

\begin{prop}
The algebras $\mathfrak{sp}(4, \mathbb{R})$ and $\mathfrak{sp}(4, \mathbb{C})$, do exist as Tractor holonomy algebras. They even exist for the `symplectic projective' construction of Section \ref{real:symplectic}.
\end{prop}
\begin{proof}
This is a sketch of a proof, without going into too many details. The Lie algebra $\mathfrak{sp}(2n+2, \mathbb{R})$ splits into
\begin{eqnarray*}
\mathfrak{sp}(2n+2, \mathbb{R}) = V^* \oplus \big( \mathfrak{sp}(2n, \mathbb{R}) \oplus \mathfrak{sp}(2, \mathbb{R}) \big) \oplus V,
\end{eqnarray*}
where $V \cong \mathbb{R}^{2n}$, $[V,V] \subset \mathfrak{sp}(2, \mathbb{R})$ and $[V^*, V^*] \subset \mathfrak{sp}(2, \mathbb{R})$. The Lie bracket between $V$ and $V^*$ is given by
\begin{eqnarray*}
\big[ X, \xi \big] = X(\xi) \cdot \nu|_{\mathfrak{sp}(2, \mathbb{R})} + X \otimes \xi + \nu(\xi) \otimes \nu(X),
\end{eqnarray*}
$\nu$ the symplectic structure. Note here that $\nu|_{\mathfrak{sp}(2, \mathbb{R})}$ is a map $V \to V^*$, equal to the identity under the isomorphism $V \cong V^*$ given by $\nu|_{\mathfrak{sp}(2n, \mathbb{R})}$, the other piece of $\nu$. We may then interpret the construction of Section \ref{real:symplectic} as a `symplectic projective structure' whose preferred connections change via
\begin{eqnarray*}
\nabla_X Y \to \nabla'_X Y = \nabla_X Y + \big[ X, \Upsilon]\cdot Y
\end{eqnarray*}
for some one-form $\Upsilon$. This implies that there exist non-flat symplectic projective manifolds in two dimensions (as $\odot^3 \mathbb{R}^2 = (\mathfrak{sp}(2,\mathbb{R}))^{(1)}$ is of dimension four, while $\mathbb{R}^{2*} = T^*_x$ is of dimension two).

Then since the tangent space of the underlying manifold $N^2$ cannot be preserved by $\overrightarrow{\nabla}$ (since $N$ cannot be Ricci-flat without being flat) we may construct a patching argument as in Proposition \ref{patch:U} to get the full tractor holonomy, using three copies patched together if need be. The process still works, as given any symplectic connection $\nabla$ with symplectic form $\nu$, and $\widehat{\nabla}$ a flat connection preserving $\nu$, then $\nabla' = f\nabla  + (1-f)\widehat{\nabla}$ also preserves $\nu$.

To generalise this argument to the complex case is slightly subtle, as we are no longer in the case of a manifold that can be made holomorphic, and the complex symplectic curvature expressions (the complex equivalent of Equations \ref{symp:curvexpr}) become considerably more complicated -- though Equations \ref{symp:curvexpr} remain valid if we look at the holomorphic ($J$-commuting) part of the curvature only.

Therefore we may start with a holomorphic symplectic connection, not $\mathbb{C}$-symplectically flat. These exist by the same argument as in the real case. Then we use partition of unity `patching' arguments on this manifold, to conjugate whatever holonomy algebra it has locally, and thus to create a manifold with full Tractor holonomy. This manifold is no longer holomorphic, but the terms from the anti-holomorphic part of the curvature cannot reduce the holonomy algebra; and since they must be contained in $\mathfrak{sp}(2n + 2, \mathbb{C})$, they can't increase it either.
\end{proof}


\begin{thebibliography}{99}
\setlength{\parskip}{0pt}
\bibitem[Ada]{E6} J.F.~Adams: \emph{Lectures on exceptional Lie groups},
Chicago Lectures in Mathematics, University of Chicago Press, Chicago, IL (1996).

\bibitem[AlMa]{quaternionic} D.V.~Alekseevsky and S.~Marchiafava: \emph{Quaternionic structures on a manifold and subordinated structures}, Ann. Mat. Pura Appl. (4) \textbf{171} (1996), 205-273.

\bibitem[Ale]{old4KM} D.V.~Alekseevsky: \emph{Riemannian spaces with unusual holonomy groups}, (Russian) Funkcional. Anal. i Prilo\v zen \textbf{2} 1968 No.~2 1-10.

\bibitem[ADM]{quaternionicweyltensor} D.V.~Alekseevsky, A.J.~Di Scala and S.~Marchiafava: \emph{Parallel K\"ahler submanifolds of quaternionic K\"ahler symmetric spaces},
http://www.hull.ac.uk/php/masdva/AToniStef051103.pdf , to be published.

\bibitem[ArLe]{me!} S.~Armstrong and T.~Leistner: \emph{Ambient connections realising conformal Tractor holonomy}, math.DG/0606410 (2006).

\bibitem[Arm1]{mecon} S.~Armstrong: \emph{Definite signature conformal holonomy: a complete classification}, math.DG/0503388 (2005).

\bibitem[Arm2]{mepro1} S.~Armstrong: \emph{Projective Holonomy I: Principles and Properties}, math.DG/0602620 (2006).

\bibitem[Arm3]{meric} S.~Armstrong: \emph{Ricci Holonomy: a Classification}, math.DG/0602619 (2006).

\bibitem[Arm4]{methesis} S.~Armstrong: \emph{Tractor Holonomy Classification for Projective and Conformal Structures}, Doctoral Thesis, Bodelian Library, Oxford University (2006).

\bibitem[Bar]{ACHHM} M.L.~Barberis: \emph{Affine connections on homogeneous hypercomplex manifolds},
J. Geom. Phys. \textbf{32} (1999), No.~1, 1-13.

\bibitem[BFGK]{SasakiEin} H.~Baum, T.~Friedrich, R.~Grunewald and I.~Kath: \emph{Twistors and Killing spinors on Riemannian manifolds},
Teubner-Texte zur Mathematik [Teubner Texts in Mathematics], \textbf{124}. B. G. Teubner Verlagsgesellschaft mbH, Stuttgart (1991).

\bibitem[BCGRS]{Symplectic connections} P.~Bieliavsky, M.~Cahen, S.~Gutt, J.~Rawnsley and L.~Schwachh\"ofer: \emph{Symplectic connections}, arXiv: math.SG/0511194v1 (2005).

\bibitem[Boh]{SasakiEin2} C.~Bohle: \emph{Killing spinors on Lorentzian manifolds}, J. Geom. Phys. \textbf{45} (2003), No.~3-4, 285-308.

\bibitem[BGN]{suexist} C.P~Boyer, K.~Galicki and M.~Nakamaye: \emph{Sasakian geometry, homotopy spheres and positive Ricci curvature}, Topology \textbf{42} (2003), No.~5, 981-1002.

\bibitem[Bry]{MEH} R.~Bryant: \emph{Metrics with exceptional holonomy},
Ann. of Math. (2) \textbf{126} (1987), No.~3, 525-576.

\bibitem[Cal]{Mob} D.M.J.~Calderbank: \emph{M\"obius structures and two-dimensional Einstein-Weyl geometry}, J. Reine Angew. Math. \textbf{504} (1998), 37-53.

\bibitem[CaGo]{ambient2} A.~\v Cap and A.R.~Gover: \emph{Standard tractors and the conformal ambient metric construction}, Ann. Global Anal. Geom. \textbf{24} (2003), No.~3, 231-259.

\bibitem[Car]{ECP} E.~Cartan: \emph{Sur les vari\'et\'es à connexion projective}, Bull. Soc. Math. France, \textbf{52} (1924), 205-241.

\bibitem[FeHi]{ambient1} C.~Fefferman and K.~Hirachi: \emph{Ambient metric construction of $Q$-curvature in conformal and CR geometries},
Math. Res. Lett. \textbf{10} (2003), No.~5-6, 819-831.

\bibitem[Fox]{fox} D.J.F.~Fox: \emph{Contact projective structures}, Indiana Univ. Math. J., \textbf{54} (2005), No.~6, 1547-1598. 

\bibitem[Joy]{CHQM} D.~Joyce: \emph{Compact hypercomplex and quaternionic manifolds}, J. Differential Geom. \textbf{35} (1992), No.~3, 743-761.

\bibitem[KoNo]{Founddiff} Kobayashi S. and Nomizu K.: \emph{Foundations of differential geometry. Vol. I}, Wiley Classics Library. A Wiley-Interscience Publication. John Wiley \& Sons, Inc., New York, 1996.

\bibitem[MeSc1]{CIH} S.~Merkulov and L.~Schwachh\"ofer: \emph{Classification of Irreducible Holonomies of Torsion-free Affine Connections}, Annals of Mathematics, \textbf{150} (1999), 77-150.

\bibitem[PPS]{HyQua} H.~Pedersen, Y.S.~Poon and A.F.~Swann: \emph{Hypercomplex structures associated to quaternionic manifolds}, Differential Geom. Appl. \textbf{9} (1998), No.~3, 273-292.

\bibitem[Sal]{DGQM} S.M.~Salamon: \emph{Differential geometry of quaternionic manifolds}, Ann. Sci. Ecole Norm. Sup. (4) \textbf{19} (1986), No.~1, 31-55.

\bibitem[Tho]{oldcone} T.Y.~Thomas: \emph{The Differential Invariants of Generalized Spaces}, National Mathematics Magazine, \textbf{9} (1935), No.~5, 151-152.

\bibitem[Wey2]{WT} H.~Weyl: \emph{Reine Infinitesimalgeometrie}, Math. Z. \textbf{2} (1918), No.~3-4, 384-411.

\end{thebibliography}
\end{document}